\documentclass[article]{IEEEtran}

\usepackage{amsmath}
\usepackage{amssymb}
\usepackage{amsfonts}
\usepackage{stfloats}
\usepackage{bbm}
\usepackage{color}
\usepackage{enumerate}	
\usepackage{wrapfig}

\usepackage{amssymb,amsfonts,amsmath,mathrsfs,mathtools,mathabx}
\usepackage{graphicx,epsfig}
\usepackage{bm}
\usepackage{multirow}
\usepackage{color}
\usepackage{xcolor}
\usepackage{subcaption}

\newtheorem{remark}{Remark}

\newcommand{\cp}{\mathbf{\bar{c}}}
\begin{document}
%
\title{Titles}
%
%
%

\title{\LARGE \bf
Modeling Large-Scale Adversarial Swarm Engagements using Optimal Control
}

\author{ 	
        Claire Walton,
         Isaac Kaminer,
         Qi Gong,        
         Abram H. Clark, and
         Theodoros Tsatsanifos

            \thanks{Claire Walton is with the Departments of Electrical Engineering and Mathematics, University of Texas at  San Antonio (e-mail:~claire.walton@utsa.edu).}
           \thanks{Isaac Kaminer and Theodoros Tsatsanifos  are with the Department of Mechanical and Aerospace Engineering,
    		 Naval Postgraduate School, Monterey, CA, 93943 USA (e-mail:~kaminer@nps.edu; ~theodoros.tsatsanifos.gr@nps.edu).}
            \thanks{Qi Gong  is with the Department of of Applied Mathematics and Statistics, University of California, Santa Cruz (e-mail:~qgong@ucsc.edu).}
    		 \thanks{Abe Clark is with the Department of Physics, Naval Postgraduate School (e-mail:~abe.clark@nps.edu).}}

\maketitle

\begin{abstract}
We study optimal control of large-scale autonomous systems under adversarial conditions in which agents may be probabilistically destroyed during an engagement. Because attrition changes the population that generates the spatial interaction dynamics, treating survival and motion independently can produce physically inconsistent solutions. We formulate a stochastic variable-population benchmark and examine three deterministic approximations that propagate agent survival probabilities from relative geometry. The models are applied to defense of a high-value unit against an attacking swarm and solved using direct optimal-control methods. Monte Carlo realizations of the stochastic benchmark show that deterministic models can reproduce engagement outcomes when attrition is coupled to spatial influence, whereas a decoupled formulation can exploit effectively destroyed defenders that continue to repel attackers, a failure mode we term ``ghost herding.'' Large-scale simulations and defender-number sweeps further show that this modeling choice can substantially alter predicted high-value-unit survival and inferred defensive resource requirements. The results provide a tractable framework for attrition-aware trajectory optimization and resource--survival analysis in adversarial autonomous systems.

\end{abstract}


%

\section{Introduction}

Recent advances in autonomous systems have enabled increasingly large teams of unmanned agents to operate cooperatively in complex environments. A central challenge is the design of coordination and control strategies that remain effective in the presence of uncertainty, disturbances, communication limitations, and loss of individual agents. These questions have motivated extensive research on robust mission planning and collective swarm behavior~\cite{evers2014robust,reynolds1987flocks,vicsek,leonard2001virtual,brambilla2013swarm}.

An additional challenge arises when autonomous agents operate in explicitly adversarial environments. In this setting, an agent's motion affects not only its ability to accomplish a mission objective, but also its exposure to an adversary and therefore its likelihood of remaining operational. Adversarial interactions may also exploit collective swarm behavior. A review of adversarial swarm control~\cite{TR2018} discusses biologically inspired strategies such as birds of prey herding bird flocks~\cite{herdingbirds}, dolphins hunting~\cite{Egerstedt}, and sheepdogs managing livestock~\cite{sheepdogs1,sheepdogs2}. Such herding strategies exploit the swarm's own response mechanisms to achieve objectives such as containment.

A related body of work has examined the nonlinear collective dynamics that arise when swarms interact or collide. Kolon and Schwartz studied collisions between self-propelled swarms and showed that communication delay, coupling strength, and collision geometry can produce qualitatively different outcomes, including flocking, milling, and scattering~\cite{Schwartz}. Hindes et al. subsequently developed an analytical approximation for colliding swarms that predicts a critical inter-swarm coupling separating scattering from the formation of a joint milling state~\cite{hindes2021critical}. Related work has also examined swarm shedding, in which coherent groups lose agents as network and dynamical conditions vary~\cite{hindes2021shedding}, while mixed-reality robotic experiments have demonstrated that theoretically predicted swarm-pattern transitions can persist on physical platforms~\cite{edwards2020delay}. Collectively, these studies show that the evolution of a multi-swarm engagement can depend sensitively on the interaction structure and on which agents participate in the collective dynamics. The present work addresses a distinct mechanism that changes this participation in adversarial engagements: stochastic attrition removes agents from the active population and therefore must also remove their subsequent influence on the spatial interaction dynamics.

More directly related to swarm defense, Chipade, Panagou, and coauthors developed the StringNet framework, in which defenders enclose and herd a risk-averse attacking swarm away from a protected region in obstacle-populated environments~\cite{chipade2019herding,chipade2021planning}. Extensions address attackers that split into multiple sub-swarms through online clustering and defender reassignment~\cite{chipade2020multiswarm,chipade2021stringnet}, three-dimensional engagements~\cite{zhang2021herding3d}, and a multi-mode strategy combining herding of risk-averse attackers with interception of risk-taking attackers~\cite{chipade2023aerial}. The StringNet approach has also been demonstrated experimentally with quadrotors~\cite{chipade2021stringnet}. These works show how defender--attacker interactions can be exploited for swarm defense. The present work focuses on an additional mechanism arising under attrition: once an agent is destroyed, it can no longer generate the interaction forces responsible for herding, avoidance, or threat. Hence survival and spatial interaction must be coupled.

A separate challenge in swarm defense is that the tactics employed by an attacking swarm may not be known a priori. Peltier et al.~\cite{peltier2025swarmclassification} used supervised neural-network time-series classification to infer adversarial swarm characteristics from observed trajectories. In that framework, two binary characteristics---inter-agent communication and the use of proportional navigation---define four distinct attacking-swarm tactics, and the classifier identifies the tactic from short histories of attacker positions and velocities. Subsequent work~\cite{peltier2025robustclassification} addressed operational uncertainty by training robust classifiers over variations in defender number, defender motion, and measurement noise, and further showed that defender trajectories can be optimized to actively elicit adversary responses that improve tactic classification. These results highlight an important sensing-and-decision aspect of swarm defense: before selecting an appropriate counter-strategy, defenders may first need to infer how the attacking swarm coordinates and responds. Rather than identifying an unknown attacker tactic, the present work examines how a specified engagement model should represent the loss of agents once attrition occurs, an issue that remains relevant regardless of how the adversarial tactic is identified.

Optimal control provides a natural framework for planning defensive trajectories in such engagements. Earlier work considered rapid-fire combat with attacker uncertainty~\cite{walton2018optimal}, while subsequent work formulated defense of a High-Value Unit (HVU) against a large-scale swarm as an uncertain-parameter optimal-control problem~\cite{walton2022parameter}. More recently, online defensive motion planning has been developed using Bernstein-polynomial approximations and model predictive control (MPC)~\cite{kang2025online}, and a distributed receding-horizon formulation has enabled individual defenders to optimize their trajectories through distributed Nash-equilibrium seeking~\cite{kang2026distributed}. These developments demonstrate the growing capability of optimization-based methods for coordinating defensive autonomous systems against adversarial swarms.

Resilience to failures and agent loss is an important topic in the broader UAV-swarm literature~\cite{phadke2022resilient}. In adversarial optimal control, however, attrition changes not only performance but also the set of agents that physically participate in the dynamics. This creates an important distinction between uncertainty in the parameters or behavior of a fixed-population dynamical system and uncertainty in the active population. If the surviving agents at time $t_k$ are represented by an index set $\mathcal{I}(t_k)$, then attrition causes $\mathcal{I}(t_k)$ to evolve stochastically. As agents are removed, the interactions among agents change and, in a reduced-state representation, so does the dimension of the dynamical system. The resulting problem therefore combines deterministic motion between attrition events with stochastic changes in the active population.

Directly incorporating this mechanism into trajectory optimization is computationally challenging, particularly for engagements involving hundreds or thousands of agents. A tempting approximation is to propagate continuous survival probabilities while retaining the original deterministic spatial dynamics. As we demonstrate in this paper, however, such a decoupling can produce qualitatively incorrect optimal solutions. In particular, an optimizer may exploit defenders having negligible survival probability as if they continued to exert repulsive influence on attacking agents. We refer to this unphysical behavior as \emph{ghost herding}.

The objective of this paper is therefore to investigate computationally tractable optimal-control models that preserve the essential coupling between attrition and spatial dynamics. To the best of the authors' knowledge, existing adversarial-swarm defense frameworks have not explicitly examined the modeling inconsistency that arises when attrition or survival is propagated while agents that have effectively been destroyed are nevertheless allowed to retain their influence on the spatial interaction dynamics. This gap is distinct from prior work on swarm herding and interception, tactic classification, parameter uncertainty, and defensive trajectory optimization. We first formulate a stochastic variable-population optimal-control problem, denoted $P_0$, in which the active-agent set changes as agents are destroyed. Because direct solution of this problem is difficult for large-scale engagements, we introduce and compare three deterministic approximations. Problem $P_1$ propagates survival probabilities but leaves the spatial dynamics uncoupled from survival; Problem $P_2$ continuously weights an agent's contribution to the interaction dynamics by its survival probability; and Problem $P_3$ removes an agent's dynamical influence once its survival probability falls below a prescribed threshold.

The principal contributions of this work are fourfold. First, we formulate an optimal-control framework in which stochastic agent attrition alters the active population and its interaction dynamics. Second, we introduce computationally tractable deterministic approximations that differ in how survival is coupled to spatial interaction. Third, we show that neglecting this coupling can generate the ghost-herding phenomenon and substantially overestimate defensive effectiveness. Finally, through large-scale numerical experiments and Monte Carlo validation, we show that the survival-dependent models $P_2$ and $P_3$ reproduce important features of the stochastic engagement much more accurately than the uncoupled model $P_1$, and we quantify how the choice of attrition model affects estimates of the defensive resources required to protect an HVU.

A preliminary version of this work appeared in our 2021 IEEE Conference on Decision and Control paper~\cite{CDC21}. That paper introduced the stochastic formulation $P_0$, the three deterministic approximations $P_1$--$P_3$, and an initial 200-realization Monte Carlo comparison. The present journal manuscript substantially extends that preliminary study. First, it develops a more explicit and internally consistent mathematical treatment of the stochastic active-agent process and of the distinction between realized survival and continuously propagated survival probability. Second, it gives explicit finite-dimensional optimal-control/NLP formulations for the deterministic approximations and clarifies how the uncoupled, survival-weighted, and thresholded models modify the spatial interactions. Third, it expands the stochastic benchmark analysis to large-scale engagements and uses common optimized defender trajectories to isolate differences among the attrition models. Fourth, the journal paper introduces a systematic force-sizing study in which defender number and relative weapon capability are varied, quantifying how the attrition model changes predicted HVU survival and the critical number of defenders required for protection. Finally, these new force-sizing results are interpreted as model-dependent resource--survival/Pareto-type frontiers, demonstrating that ghost herding can produce an artificially optimistic resource--performance trade-off. The journal contribution is therefore an expanded mathematical, computational, and resource-allocation analysis of the preliminary framework in~\cite{CDC21}, rather than a republication of the conference results.

\section{Modeling and Optimization Framework for Adversarial Swarming} \label{sec:sec_3.1}
\label{sec:modeling}

We consider multi-agent systems for which the motion of each active agent is deterministic between attrition events, whereas the loss of agents is stochastic. Thus, conditioned on the set of agents that are operational at a given time, the local motion model contains no process noise; randomness enters through the attrition process. This distinction is important because the loss of an agent changes not only a survival metric but also the interaction structure of the remaining swarm.

Let $\mathcal I(t_k)$ denote the set of active agents at time $t_k$. For an engagement with $N$ attackers and $M$ defenders, $\mathcal I(t_0)$ contains all $N+M$ agents, and $\mathcal I(t_{k+1})\subseteq\mathcal I(t_k)$ because agents that are lost do not re-enter the engagement. Let
\[
z(t_k)=\operatorname{col}_{i\in\mathcal I(t_k)} z_i(t_k),\qquad
u(t_k)=\operatorname{col}_{i\in\mathcal I(t_k)} u_i(t_k),
\]
where $z_i$ is the state of agent $i$ and $u_i$ is its control input, when applicable. Conditioned on $\mathcal I(t_k)$, the discrete-time dynamics can be written as
\begin{equation}
 z(t_{k+1}) = \phi^k\!\left(z(t_k),u(t_k),\mathcal I(t_k)\right).
\label{eq:agent_dyn}
\end{equation}
The superscript $k$ emphasizes that the collective dynamics change when the active-agent set changes. In a reduced-state representation, the dimension of $z(t_k)$ is therefore random. Equivalently, one may retain a fixed-dimensional state and multiply each agent's contribution by a binary alive/dead indicator; both descriptions represent the same variable-population mechanism.

To distinguish a survival probability from an actual survival realization, define
\[
p_j(t_k)\in[0,1]
\]
as the \emph{conditional one-step survival probability} of active agent $j$ on $[t_k,t_{k+1}]$, given the current state and active-agent set. In general,
\begin{equation}
 p(t_k)=\mathcal P\!\left(z(t_k),\mathcal I(t_k)\right),
\label{eq:conditional_survival}
\end{equation}
where $p(t_k)$ collects the one-step survival probabilities of the currently active agents. Let $\omega_j(t_k)\sim\mathcal U[0,1]$ be independent random variables and define the attrition set
\begin{equation}\label{J eqn}
 J(t_k)\triangleq
 \left\{j\in\mathcal I(t_k)\;\middle|\;\omega_j(t_k)>p_j(t_k)\right\}.
\end{equation}
The active-agent set then evolves according to
\begin{equation}\label{update eqn}
 \mathcal I(t_{k+1})
 =\psi\!\left(\mathcal I(t_k),z(t_k),\omega(t_k)\right)
 \triangleq \mathcal I(t_k)\setminus J(t_k).
\end{equation}
Equivalently,
\begin{equation}
\mathcal I(t_{k+1})=\psi\!\left(\mathcal I(t_k),p(t_k),\omega(t_k)\right),
\label{eq:I_dyn}
\end{equation}
with $p(t_k)$ determined by \eqref{eq:conditional_survival}. This formulation makes explicit the distinction between the probability $p_j(t_k)$ and the binary event of whether agent $j$ actually survives the interval.

Because the active-agent set is random, both the terminal state and any terminal performance measure are random. A general stochastic optimal-control problem for adversarial swarming can therefore be written as
\begin{equation}
P_0\triangleq
\left\{
\begin{array}{ll}
\displaystyle \min_{u} & J:=\mathbb E\!\left[F\!\left(z(t_f),\mathcal I(t_f)\right)\right] \\
\textbf{subject to}
& z(t_{k+1})=\phi^k\!\left(z(t_k),u(t_k),\mathcal I(t_k)\right),\\
& p(t_k)=\mathcal P\!\left(z(t_k),\mathcal I(t_k)\right),\\
& \mathcal I(t_{k+1})=\psi\!\left(\mathcal I(t_k),p(t_k),\omega(t_k)\right),\\
& \mathcal H\!\left(z(t_k),u(t_k)\right)\le 0.
\end{array}
\right.
\label{OCP_A}
\end{equation}
Here $F$ is the terminal cost, $\phi^k$ denotes the deterministic motion of the active population, $\mathcal P$ determines conditional survival probabilities, and $\mathcal H$ collects state and control constraints.

Problem $P_0$ is difficult for three related reasons. First, attrition causes stochastic switching of the active population and, in a reduced-state description, a random time-varying state dimension. Second, the attrition process is coupled to the motion because survival probabilities depend on relative geometry, while the subsequent motion depends on which agents survived. Third, large engagements lead to high-dimensional trajectory-optimization problems even before stochastic branching is considered. The purpose of the approximations introduced next is therefore not merely to replace randomness by a deterministic survival metric, but to retain--at least approximately--the essential feedback from attrition to spatial interaction.

\section{Proposed Solution Methods} \label{sec:soving-zombie}

In this section, we introduce three deterministic approximations to Problem $P_0$. The key distinction among them is how the effect of attrition is fed back into the spatial dynamics.

\subsection{$P_1$: Deterministic and Decoupled Optimal Control Problem}

We first consider the approximation used in~\cite{walton2018optimal}, in which all agents remain present in the motion model throughout the engagement. Attrition is represented by deterministic survival-probability states rather than by a random active-agent set. For the HVU-defense example developed in the next section, let $Q_j^A(t)$ denote the modeled probability that attacker $j$ survives to time $t$, $Q_l^D(t)$ the corresponding probability for defender $l$, and $Q_0(t)$ the probability that the HVU survives. Collect these quantities as
\[
Q(t)=\left[Q_0(t),Q_1^A(t),\ldots,Q_N^A(t),Q_1^D(t),\ldots,Q_M^D(t)\right]^T.
\]
Their deterministic evolution is written abstractly as
\begin{equation}
Q(t_{k+1})=\Psi\!\left(Q(t_k),z(t_k)\right).
\label{eq:Q_abstract}
\end{equation}
Unlike the random set $\mathcal I(t_k)$ in $P_0$, $Q(t_k)$ is a deterministic continuous-valued surrogate once the trajectories are specified.

The resulting optimal-control problem is
\begin{equation} \label{OCP}
P_1\triangleq
\left\{
\begin{array}{ll}
\displaystyle \min_u & J:=F\!\left(z(t_f),Q(t_f)\right)\\
\textbf{subject to}
& z(t_{k+1})=\phi\!\left(z(t_k),u(t_k)\right),\\
& Q(t_{k+1})=\Psi\!\left(Q(t_k),z(t_k)\right),\\
& \mathcal H\!\left(z(t_k),u(t_k)\right)\le 0.
\end{array}
\right.
\end{equation}
Problem $P_1$ is computationally convenient because its state dimension is fixed and its dynamics do not switch when survival probabilities become small. That same decoupling is also its principal weakness: an agent with negligible survival probability can continue to exert its full dynamical influence. In the HVU-defense example, the optimizer can exploit this inconsistency and produce ``ghost-herding'' solutions in which effectively destroyed defenders continue to repel attackers.

\subsection{$P_2$: Survival-Weighted Interaction Model}

The second approximation feeds survival back into the motion model continuously. Each agent's contribution to the interaction dynamics is weighted by its modeled survival probability. In abstract form,
\begin{equation} \label{P2}
P_2\triangleq
\left\{
\begin{array}{ll}
\displaystyle \min_u & J:=F\!\left(z(t_f),Q(t_f)\right)\\
\textbf{subject to}
& z(t_{k+1})=\phi\!\left(z(t_k),u(t_k),Q(t_k)\right),\\
& Q(t_{k+1})=\Psi\!\left(Q(t_k),z(t_k)\right),\\
& \mathcal H\!\left(z(t_k),u(t_k)\right)\le 0.
\end{array}
\right.
\end{equation}
If all survival probabilities are equal to one, $P_2$ reduces to $P_1$. If the weights were binary alive/dead indicators, the interaction structure would coincide with that of a particular realization of $P_0$. Thus, $P_2$ can be interpreted as a smooth relaxation of binary participation. It is not identical to the stochastic model--for example, attraction between two attackers is gradually weakened as their survival probabilities decrease--but it removes the central ghost-herding inconsistency because an agent's dynamical influence vanishes as its survival probability approaches zero.

\subsection{$P_3$: Thresholded Interaction Model}

The third approximation replaces the random active-agent set by a deterministic thresholded set. Let $\mathcal I^m(t_k)$ denote the agents whose modeled survival probability remains above a prescribed threshold $q_{\rm th}\in(0,1)$; in the numerical study we use $q_{\rm th}=0.5$. Agents in $\mathcal I^m(t_k)$ interact at full strength, while agents outside this set no longer contribute to either spatial interactions or mutual attrition. The resulting problem is
\begin{equation} \label{P3}
P_3\triangleq
\left\{
\begin{array}{ll}
\displaystyle \min_u & J:=F\!\left(z(t_f),Q(t_f),\mathcal I^m(t_f)\right)\\
\textbf{subject to}
& z(t_{k+1})=\phi^k\!\left(z(t_k),u(t_k),\mathcal I^m(t_k)\right),\\
& Q(t_{k+1})=\Psi\!\left(Q(t_k),z(t_k),\mathcal I^m(t_k)\right),\\
& \mathcal I^m(t_{k+1})=\psi_m\!\left(\mathcal I^m(t_k),Q(t_{k+1})\right),\\
& \mathcal H\!\left(z(t_k),u(t_k)\right)\le 0.
\end{array}
\right.
\end{equation}
The use of $Q(t_{k+1})$ in the threshold update means that the survival state is first propagated over $[t_k,t_{k+1}]$ and agents crossing the threshold are removed before the next interval. Unlike $P_0$, the update is deterministic; unlike $P_2$, participation is binary rather than continuously weighted.

\subsection{Validation by Monte Carlo Simulation}

To evaluate the deterministic approximations, we use the optimized defender trajectories as fixed inputs to an ensemble of Monte Carlo simulations of the stochastic attrition process. In each realization, spatial dynamics are propagated only for active agents, and an agent that is destroyed no longer contributes to subsequent interactions or attrition.

For active agent $j$, let $p_j(t_k)$ denote its conditional probability of surviving the interval $[t_k,t_{k+1}]$, computed from the current geometry and the currently active opponents. Draw independent $\omega_j(t_k)\sim\mathcal U[0,1]$ and define
\begin{equation}\label{eq:J_MC}
J_{\rm MC}(t_k)\triangleq
\left\{j\in\mathcal I(t_k)\;\middle|\;\omega_j(t_k)>p_j(t_k)\right\}.
\end{equation}
The active set is then updated by
\[
\mathcal I(t_{k+1})=\mathcal I(t_k)\setminus J_{\rm MC}(t_k).
\]
When a cumulative survival variable $Q_j$ is propagated using the active-agent realization, the same one-step probability may be written as
\[
p_j(t_k)=\frac{Q_j(t_{k+1})}{Q_j(t_k)},
\]
provided $Q_j(t_k)>0$. This ratio is a \emph{conditional one-step} probability and should not be confused with the binary survival realization itself.

The Monte Carlo ensemble therefore provides a direct numerical benchmark for the attrition--dynamics coupling represented by $P_0$, while avoiding the exponential scenario branching that would arise if the random attrition events were included directly inside a gradient-based trajectory optimization.

\section{Case Study: Defense Against a Swarm Attack}\label{sec:model}

To evaluate the models introduced above, we consider a swarm attacking a high-value unit (HVU). The HVU is protected by $M$ defenders whose trajectories are optimized to maximize the HVU survival probability. For the deterministic approximations $P_1$--$P_3$, the terminal cost is
\[
F(z(t_f),Q(t_f))=1-Q_0(t_f),
\]
where $Q_0(t_f)$ is the modeled probability that the HVU survives to the final time. For the stochastic benchmark $P_0$, the corresponding objective is the probability of HVU destruction, $1-\mathbb E[\chi_0(t_f)]$, where $\chi_0(t)\in\{0,1\}$ is the HVU alive indicator.

\subsection{Attacker Equations of Motion}

We consider $N$ attacking agents and $M$ defending agents. Attacker $i$ has position $x_i(t)\in\mathbb R^3$ and velocity $v_i(t)\in\mathbb R^3$, while defender $l$ has position $s_l(t)\in\mathbb R^3$ and velocity $v_{D_l}(t)\in\mathbb R^3$. In the uncoupled model, the continuous-time dynamics of attacker $i$ are
\begin{align}\label{eqn:attackers}
\dot x_i &= v_i,\nonumber\\
\dot v_i &= \sum_{\substack{j=1\\j\neq i}}^N
\frac{f_I(x_{ij})}{\|x_{ij}\|}x_{ij}
+\sum_{l=1}^M \frac{f_d(s_{il})}{\|s_{il}\|}s_{il}
+K\frac{h_i}{\|h_i\|}-b v_i
=: f_i(t),
\end{align}
where $x_{ij}=x_i-x_j$, $s_{il}=x_i-s_l$, $h_i=h-x_i$, and $h$ is the HVU position. The four terms represent attacker--attacker interaction, repulsion from defenders, attraction toward the HVU through a virtual leader, and velocity damping, respectively.

For $f_I$ and $f_d$ we use the model of Leonard and Fiorelli~\cite{leonard2001virtual}. The attacker--attacker interaction is repulsive for $\|x_{ij}\|\le d_0$, attractive for $d_0<\|x_{ij}\|\le d_1$, and zero for $\|x_{ij}\|>d_1$. The attacker--defender interaction is purely repulsive for $\|s_{il}\|\le s_0$ and zero beyond that range. Simulations using the Reynolds model~\cite{reynolds1987flocks} in place of the first interaction term produced qualitatively similar behavior.

\subsection{Defender Equations of Motion}

Each defender is modeled as a three-dimensional double integrator,
\begin{align}\label{eqn:defenders}
\dot s_l &= v_{D_l},\nonumber\\
\dot v_{D_l} &= u_l,
\end{align}
where $u_l(t)\in\mathbb R^3$ and $\|u_l(t)\|_\infty\le u_{\max}$.

\subsection{Mutual Attrition Model}

\begin{figure}
\raggedright
(a) \\ \centering \includegraphics[trim=5mm 2mm 5mm 5mm, clip, width=0.8\columnwidth]{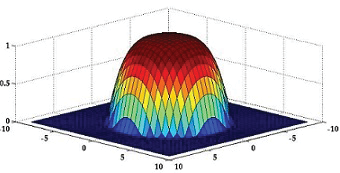} \\
\caption{Normal-distribution-based function $\Phi$ used to define the distance-dependent attrition rates.}
\label{fig:poission}
\end{figure}

Let $d_{il}^{\rm att}$ denote the attrition rate of attacker $i$ due to defender $l$, $d_{li}^{\rm def}$ the attrition rate of defender $l$ due to attacker $i$, and $d_i^{\rm hvu}$ the attrition rate of the HVU due to attacker $i$. We use the distance-dependent functions
\begin{align}
 d_{il}^{\rm att} &= \lambda_d\Phi\!\left(\frac{\|s_{il}\|^2}{\sigma_d}\right),\nonumber\\
 d_{li}^{\rm def} &= \lambda_a\Phi\!\left(\frac{\|s_{il}\|^2}{\sigma_a}\right),\nonumber\\
 d_i^{\rm hvu} &= \lambda_a\Phi\!\left(\frac{\|h_i\|^2}{\sigma_a}\right),
\label{eq:damage_rates}
\end{align}
where $\sigma_a,\sigma_d$ characterize weapon range and $\lambda_a,\lambda_d$ characterize firing rate. The function $\Phi$ is shown in Fig.~\ref{fig:poission}; it is largest at zero separation and decreases smoothly with distance.

For the deterministic models, the cumulative survival probabilities satisfy
\begin{align}
Q_i^A(t_{k+1})
&=Q_i^A(t_k)\prod_{l=1}^M
\left(1-d_{il}^{\rm att}(t_k)Q_l^D(t_k)\Delta t\right),\label{eq:Q_i(t)}\\
Q_l^D(t_{k+1})
&=Q_l^D(t_k)\prod_{i=1}^N
\left(1-d_{li}^{\rm def}(t_k)Q_i^A(t_k)\Delta t\right),\nonumber\\
Q_0(t_{k+1})
&=Q_0(t_k)\prod_{i=1}^N
\left(1-d_i^{\rm hvu}(t_k)Q_i^A(t_k)\Delta t\right).
\label{eqn:Q(t)-all}
\end{align}
All survival probabilities are initialized at one. These equations make the indexing explicit: attacker survival depends on defenders, defender survival depends on attackers, and HVU survival depends on the attackers' HVU damage rates.

\subsection{Problem Formulation}

Let $t_{k+1}=t_k+\Delta t$, $k=0,\ldots,K-1$, with $t_0=0$ and $t_K=t_f$. We now specialize $P_0$--$P_3$ to the HVU-defense problem.

\subsubsection{Stochastic benchmark $P_0$}

Let $\mathcal I_A(t_k)\subseteq\{1,\ldots,N\}$ and $\mathcal I_D(t_k)\subseteq\{1,\ldots,M\}$ denote the active attackers and defenders. The attacker acceleration for a surviving attacker $i\in\mathcal I_A(t_k)$ is
\begin{align}
f_i^{\mathcal I}(t_k)=
&\sum_{\substack{j\in\mathcal I_A(t_k)\\j\neq i}}
\frac{f_I(x_{ij})}{\|x_{ij}\|}x_{ij}
+\sum_{l\in\mathcal I_D(t_k)}
\frac{f_d(s_{il})}{\|s_{il}\|}s_{il}\nonumber\\
&+K\frac{h_i}{\|h_i\|}-b v_i.
\label{eq:f_active}
\end{align}
Conditioned on the current realization, the one-step survival probabilities are
\begin{align}
p_i^A(t_k)&=\prod_{l\in\mathcal I_D(t_k)}
\left(1-d_{il}^{\rm att}(t_k)\Delta t\right),\nonumber\\
p_l^D(t_k)&=\prod_{i\in\mathcal I_A(t_k)}
\left(1-d_{li}^{\rm def}(t_k)\Delta t\right),\nonumber\\
p_0(t_k)&=\prod_{i\in\mathcal I_A(t_k)}
\left(1-d_i^{\rm hvu}(t_k)\Delta t\right).
\label{eq:p0_conditional}
\end{align}
Independent uniform random variables are then used to determine the binary survival realization over the interval, as described in Section~\ref{sec:soving-zombie}. Thus, $P_0$ minimizes
\begin{equation}
J_0=1-\mathbb E[\chi_0(t_f)]
\label{eq:P0_cost}
\end{equation}
subject to the active-agent dynamics, the defender control constraints, and the stochastic active-set updates. This is the benchmark model represented by the Monte Carlo simulations.

\subsubsection{Decoupled deterministic model $P_1$}

In $P_1$, all $N$ attackers and $M$ defenders remain in the spatial dynamics for the entire horizon. The survival states are propagated using \eqref{eq:Q_i(t)}--\eqref{eqn:Q(t)-all}, and the objective is
\[
J_1=1-Q_0(t_f).
\]
The optimization is over the defender controls $u_l(t)$, $l=1,\ldots,M$, subject to \eqref{eqn:attackers}, \eqref{eqn:defenders}, and $\|u_l(t)\|_\infty\le u_{\max}$. Because the dynamics do not depend on $Q$, even a defender with $Q_l^D\approx0$ continues to repel attackers.

\subsubsection{Survival-weighted model $P_2$}

In $P_2$, the attacker acceleration is replaced by
\begin{align}\label{eqn:wieighted}
f_i^W(t_k)=
&\sum_{\substack{j=1\\j\neq i}}^N Q_j^A(t_k)
\frac{f_I(x_{ij})}{\|x_{ij}\|}x_{ij}
+\sum_{l=1}^M Q_l^D(t_k)
\frac{f_d(s_{il})}{\|s_{il}\|}s_{il}\nonumber\\
&+K\frac{h_i}{\|h_i\|}-b v_i,
\qquad i=1,\ldots,N.
\end{align}
The survival probabilities are still propagated by \eqref{eq:Q_i(t)}--\eqref{eqn:Q(t)-all}. Hence $P_2$ uses a continuous relaxation of participation: the effect of agent $j$ on the motion model decreases continuously with its survival probability.

\subsubsection{Threshold model $P_3$}

For $P_3$, define binary participation weights
\begin{align}
W_i^A(t_k)&=\mathbbm 1\!\left\{Q_i^A(t_k)>q_{\rm th}\right\},\nonumber\\
W_l^D(t_k)&=\mathbbm 1\!\left\{Q_l^D(t_k)>q_{\rm th}\right\},
\label{eq:threshold_weights}
\end{align}
with $q_{\rm th}=0.5$. Equivalently, these weights define deterministic active sets $\mathcal I_A^m(t_k)$ and $\mathcal I_D^m(t_k)$. The attacker acceleration becomes
\begin{align}\label{eqn:threshold}
f_i^m(t_k)=
&\sum_{\substack{j=1\\j\neq i}}^N W_j^A(t_k)
\frac{f_I(x_{ij})}{\|x_{ij}\|}x_{ij}
+\sum_{l=1}^M W_l^D(t_k)
\frac{f_d(s_{il})}{\|s_{il}\|}s_{il}\nonumber\\
&+K\frac{h_i}{\|h_i\|}-b v_i.
\end{align}
The survival-probability updates are modified consistently so that only participating opponents can inflict attrition:
\begin{align}
Q_i^A(t_{k+1})
&=Q_i^A(t_k)\prod_{l=1}^M
\left(1-d_{il}^{\rm att}(t_k)W_l^D(t_k)\Delta t\right),\nonumber\\
Q_l^D(t_{k+1})
&=Q_l^D(t_k)\prod_{i=1}^N
\left(1-d_{li}^{\rm def}(t_k)W_i^A(t_k)\Delta t\right),\nonumber\\
Q_0(t_{k+1})
&=Q_0(t_k)\prod_{i=1}^N
\left(1-d_i^{\rm hvu}(t_k)W_i^A(t_k)\Delta t\right).
\label{eq:P3_survival}
\end{align}
For compatibility with the earlier notation, the set of agents removed at step $k$ may be written as
\begin{equation}\label{JT eqn}
J^m(t_k)=\left\{j:\;Q_j(t_k)\le q_{\rm th}\right\}.
\end{equation}
Unlike the random attrition set in $P_0$, $J^m(t_k)$ is completely determined by the propagated survival probabilities.

The numerical transcription and trajectory parameterization used to solve $P_1$--$P_3$ are described in the next section.

\subsection{Numerical Methods}

We use direct methods of optimal control~\cite{cichellaBernsteinCDC,cichellaBernsteinIEEETAC} to obtain numerical solutions to the deterministic problems $P_1$, $P_2$, and $P_3$. The optimization variables parameterize the defender trajectories, and the terminal objective is the probability of HVU destruction, $1-Q_0(t_f)$. To scale the forward simulation to large numbers of attackers and defenders, we use numerical techniques motivated by molecular dynamics~\cite{allen2017computer}. In particular, the attacker equations of motion are integrated using the same velocity-Verlet-based implementation used in our previous simulations~\cite{Verlet1967,plimpton1995fast}. The attrition states are propagated with the model-specific update equations defined in Section~\ref{sec:model}.

Let $s_{m_L}$ denote an $L$th-order Bernstein-polynomial approximation of defender position $s_m$.  With control points $\cp_{mj}$ and Bernstein basis $b_{j,L}(t)$ on $[0,t_f]$, the position, its $r$th derivative, and the double-integrator control are
\begin{align}
 s_{m_L}(t)&=\sum_{j=0}^{L}\cp_{mj}b_{j,L}(t), \label{eq:beziercurve}\\
 b_{j,L}(t)&=\binom{L}{j}\frac{t^j(t_f-t)^{L-j}}{t_f^L}, \nonumber\\
 s^{(r)}_{m_L}(t)&=\sum_{j=0}^{L}\left(\sum_{i=0}^{L}\cp_{mi}D^r_{ij}\right)b_{j,L}(t), \label{eq:derivatives}\\
 u_{m_L}(t)&=\ddot s_{m_L}(t)
 =\sum_{j=0}^{L}\left(\sum_{i=0}^{L}\cp_{mi}D^2_{ij}\right)b_{j,L}(t). \label{eq:2derivative}
\end{align}
Here $D^r_{ij}$ is the $(i,j)$ entry of the $r$th power of the Bernstein differentiation matrix $\mathbf D\in\mathbb R^{(L+1)\times(L+1)}$; see~\cite[Chapter 5]{farouki2012bernstein}.  Parameterizing position directly is sufficient because the defender dynamics are double integrators.

Using the Bernstein parameterization above, the deterministic optimal-control problem becomes a finite-dimensional nonlinear program in the control points.  For $P_1$, we use the standard velocity--Verlet position and velocity recursions for the attacker states.  Writing $x_{i,k}=x_i(t_k)$, $v_{i,k}=v_i(t_k)$, $Q^A_{i,k}=Q_i^A(t_k)$, $Q^D_{m,k}=Q_m^D(t_k)$, and $Q_{0,k}=Q_0(t_k)$, the NLP is
\begingroup
\small
\begin{equation}
\label{eq:P1_NLP}
\begin{aligned}
P_1^L:\qquad
\min_{\{\cp_{m\ell}\}}\quad
& J^L=1-Q_{0,K} \\
\text{s.t.}\quad
&s_{m,k}=\sum_{\ell=0}^{L}\cp_{m\ell}b_{\ell,L}(t_k), \\
&u_{m,k}=\sum_{j=0}^{L}
\left(\sum_{\ell=0}^{L}\cp_{m\ell}D^2_{\ell j}\right)b_{j,L}(t_k), \\
&x_{i,k+1}=x_{i,k}+v_{i,k}\Delta t
+\tfrac{1}{2}f_{i,k}\Delta t^2, \\
&v_{i,k+1}=v_{i,k}
+\tfrac{1}{2}\bigl(f_{i,k}+f_{i,k+1}\bigr)\Delta t, \\
&f_{i,k}=\sum_{\substack{j=1\\j\neq i}}^{N}
\frac{f_I(x_{ij,k})}{\|x_{ij,k}\|}x_{ij,k} \\
&\qquad+\sum_{m=1}^{M}
\frac{f_d(s_{im,k})}{\|s_{im,k}\|}s_{im,k} \\
&\qquad+K\frac{h_{i,k}}{\|h_{i,k}\|}-b v_{i,k}, \\
&Q^A_{i,k+1}=Q^A_{i,k}
\prod_{m=1}^{M}\!\left(1-d_{im,k}^{\rm att}Q^D_{m,k}\Delta t\right), \\
&Q^D_{m,k+1}=Q^D_{m,k}
\prod_{i=1}^{N}\!\left(1-d_{mi,k}^{\rm def}Q^A_{i,k}\Delta t\right), \\
&Q_{0,k+1}=Q_{0,k}
\prod_{i=1}^{N}\!\left(1-d_{i,k}^{\rm hvu}Q^A_{i,k}\Delta t\right), \\
&\|u_{m,k}\|_\infty\le u_{\max}, \\
&\|s_{m,k}-s_{\ell,k}\|\ge d_{\min},\quad m\neq \ell.
\end{aligned}
\end{equation}
\endgroup
where $i=1,\ldots,N$, $m,\ell=1,\ldots,M$, and $k=0,\ldots,K-1$, with the prescribed initial attacker states and $Q^A_{i,0}=Q^D_{m,0}=Q_{0,0}=1$.  The quantities $x_{ij,k}=x_{i,k}-x_{j,k}$, $s_{im,k}=x_{i,k}-s_{m,k}$, and $h_{i,k}=h-x_{i,k}$ are evaluated at the corresponding grid point.  Thus, for a candidate set of Bernstein control points, all attacker and survival states are obtained by forward propagation and the terminal cost is evaluated directly.

For completeness, the finite-dimensional formulations for $P_2$ and $P_3$ use the same objective, Bernstein decision variables, state-update equations, control bounds, and separation constraints, but with the following explicit model substitutions.  For $P_2$, the acceleration constraint in \eqref{eq:P1_NLP} is replaced by
\begin{align}
 f^W_{i,k}={}&\sum_{\substack{j=1\\j\neq i}}^{N}Q^A_{j,k}
 \frac{f_I(x_{ij,k})}{\|x_{ij,k}\|}x_{ij,k}\nonumber\\
 &+\sum_{m=1}^{M}Q^D_{m,k}
 \frac{f_d(s_{im,k})}{\|s_{im,k}\|}s_{im,k}\nonumber\\
 &+K\frac{h_{i,k}}{\|h_{i,k}\|}-b v_{i,k},
\label{eq:P2_NLP_force}
\end{align}
and $f_{i,k}$ and $f_{i,k+1}$ in the velocity--Verlet recursion are replaced by $f^W_{i,k}$ and $f^W_{i,k+1}$, respectively.  The survival recursions remain the three equations shown explicitly in \eqref{eq:P1_NLP}.

For $P_3$, define
\[
W^A_{i,k}=\mathbbm 1\{Q^A_{i,k}>q_{\rm th}\},\qquad
W^D_{m,k}=\mathbbm 1\{Q^D_{m,k}>q_{\rm th}\}.
\]
The acceleration constraint becomes
\begin{align}
 f^m_{i,k}={}&\sum_{\substack{j=1\\j\neq i}}^{N}W^A_{j,k}
 \frac{f_I(x_{ij,k})}{\|x_{ij,k}\|}x_{ij,k}\nonumber\\
 &+\sum_{m=1}^{M}W^D_{m,k}
 \frac{f_d(s_{im,k})}{\|s_{im,k}\|}s_{im,k}\nonumber\\
 &+K\frac{h_{i,k}}{\|h_{i,k}\|}-b v_{i,k},
\label{eq:P3_NLP_force}
\end{align}
while the survival recursions are
\begin{align}
Q^A_{i,k+1}&=Q^A_{i,k}\prod_{m=1}^{M}
\left(1-d_{im,k}^{\rm att}W^D_{m,k}\Delta t\right),\nonumber\\
Q^D_{m,k+1}&=Q^D_{m,k}\prod_{i=1}^{N}
\left(1-d_{mi,k}^{\rm def}W^A_{i,k}\Delta t\right),\nonumber\\
Q_{0,k+1}&=Q_{0,k}\prod_{i=1}^{N}
\left(1-d_{i,k}^{\rm hvu}W^A_{i,k}\Delta t\right).
\label{eq:P3_NLP_survival}
\end{align}
Accordingly, $P_1^L$, $P_2^L$, and $P_3^L$ have the same finite-dimensional decision vector and geometric constraints; they differ only in how the survival state modifies the spatial interaction and attrition terms.

\begin{remark}
The double-integrator defender model is differentially flat with flat output $s_m$, so $u_m(t)=\ddot s_m(t)$. Parameterizing the defender position directly therefore eliminates the need to include separate defender-state propagation constraints in the NLP.
\end{remark}

\begin{figure}
\raggedright 
(a) \\ \centering \includegraphics[trim=5mm 2mm 5mm 5mm, clip, width=0.8\columnwidth]{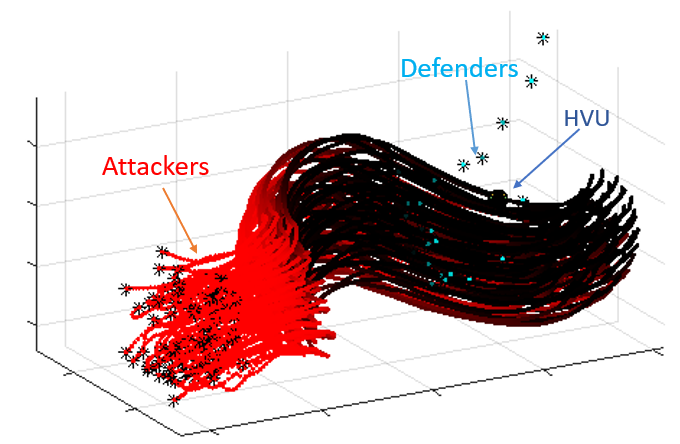} \\
\raggedright (b) \\ \centering \includegraphics[trim=2mm 2mm 5mm 5mm, clip, width=0.8\columnwidth]{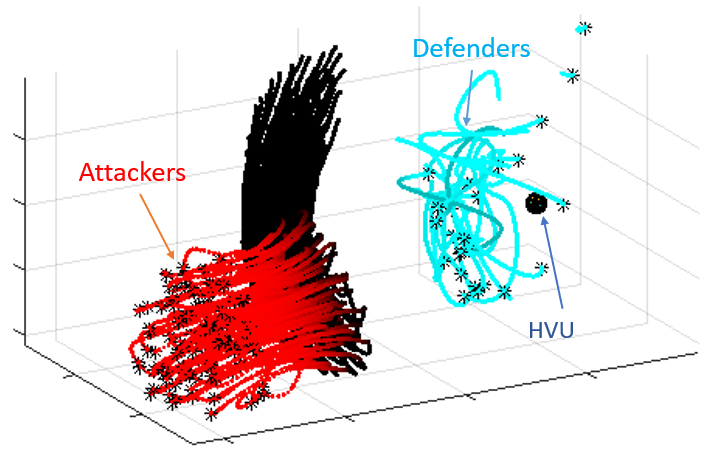} \\
\raggedright(c) \\ \centering \includegraphics[width=\columnwidth]{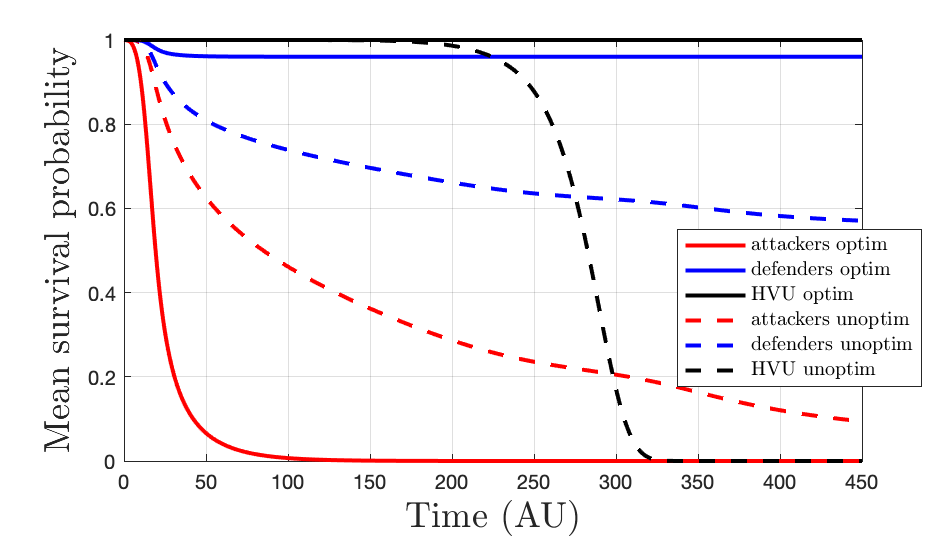}
\caption{Unoptimized (a) and optimized (b) defender trajectories for 100 attacking swarm agents and 25 defenders with superior weapons; the optimized trajectories improve HVU defense, as shown in (c). Reproduced from~\cite{CDC21}. \textcopyright~2021 IEEE.}
\label{fig:opt-unopt}
\end{figure}

Figure~\ref{fig:opt-unopt} shows an optimization of $P_1$, where survival probability is not coupled to the equations of motion of the attackers, for 25 defenders protecting an HVU against a swarm of 100 attackers. The defenders have a 50\% larger weapons range ($\sigma_d/\sigma_a = 1.5$) as well as double the fire rate with respect to the attackers ($\lambda_d/\lambda_a = 2$). Figure~\ref{fig:opt-unopt}(a) shows the results of a simulation where the defenders remain in place, and Figure~\ref{fig:opt-unopt}(b) shows results of a simulation after optimizing the trajectories of the defenders. The attackers are red and the defenders are cyan, but turn to black as their survival probability decreases to zero. In both scenarios the defenders are stronger and suffer fewer losses; however, in the unoptimized scenario, some of the attacking agents manage to penetrate the defenders’ zone and destroy the HVU. Figure \ref{fig:opt-unopt}(c) shows that the optimized trajectories lead to better results for the survival of the HVU, the survival of defenders, and the defeat of the attackers (although the survival of the HVU is the only metric used in the optimization). We also note that the defending agents utilize both herding and weapons, as seen in Fig.~\ref{fig:opt-unopt}.

Similar results to those shown in Fig.~\ref{fig:opt-unopt} are found when considering problems $P_2$ (where spatial interactions are multiplied by the survival probabilities) and $P_3$ (where spatial interactions and damage functions are turned off below a 50\% survival probability threshold). These results are qualitatively indistinguishable from the results we show for $P_1$, so we do not show them here. However, a key question remains: how well do the modeling frameworks used in problems $P_1$ (uncoupled dynamics), $P_2$, and $P_3$ compare with the stochastic problem $P_0$? Problems $P_1$, $P_2$, and $P_3$ are approximations to $P_0$. 


\subsection{Comparing  Performance of the Proposed Models}\label{sec:optimization comparizon}



To answer the question posed in the previous section, we compare numerical solutions of problems $P1 - P3$ with the outcome of a Monte Carlo simulation with fixed defender trajectories  We consider a case with $M=200$ defenders and $N=2066$ attackers with identical rates and ranges of fire ($\lambda_a = \lambda_d$,  $\sigma_a = \sigma_d$). The defender trajectories are obtained by solving problems $P1,P2,P3$ and these trajectories are then tested using Monte Carlo simulation. Results for Monte Carlo simulation. are averaged over 200 runs; all other results use a single simulation, since there is no randomness.

\begin{figure*}
\centering
\includegraphics[trim=0mm 0mm 0mm 0mm, clip, width=0.4\textwidth]{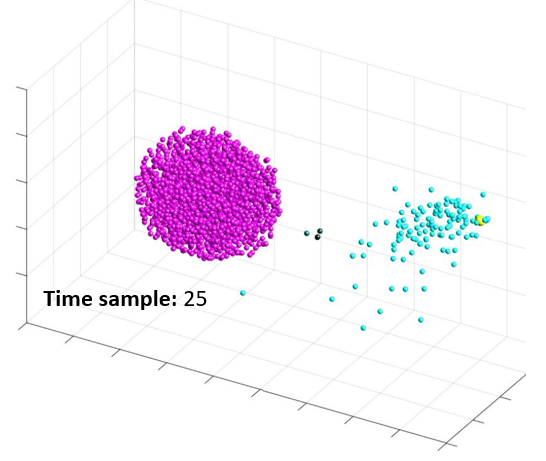}
\includegraphics[trim=0mm 0mm 0mm 0mm, clip,width=0.32\textwidth]{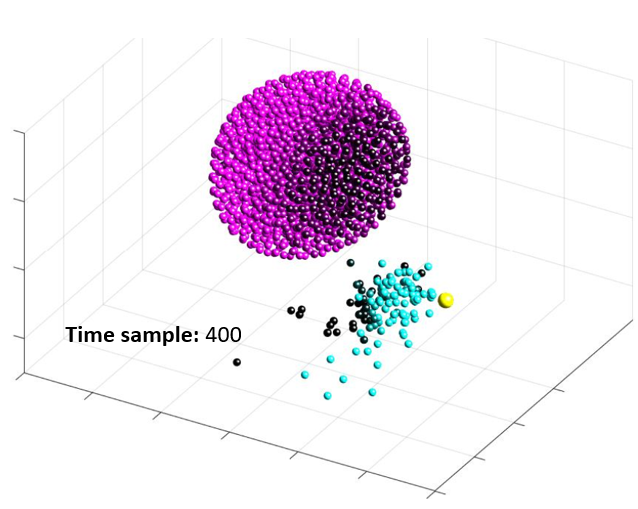}
\includegraphics[trim=0mm 0mm 0mm 0mm, clip,width=0.225\textwidth]{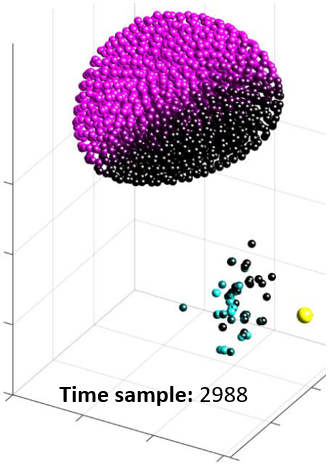}
\includegraphics[width=0.33\textwidth]{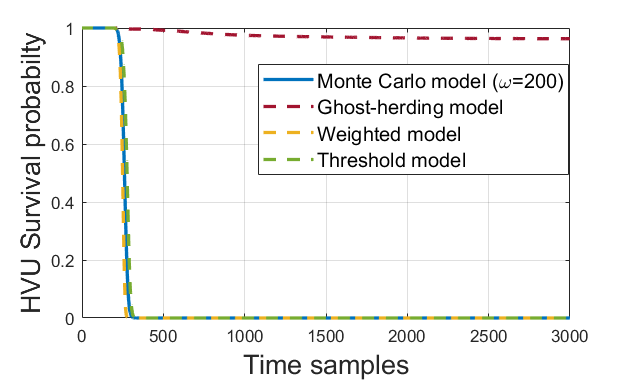}
\includegraphics[width=0.33\textwidth]{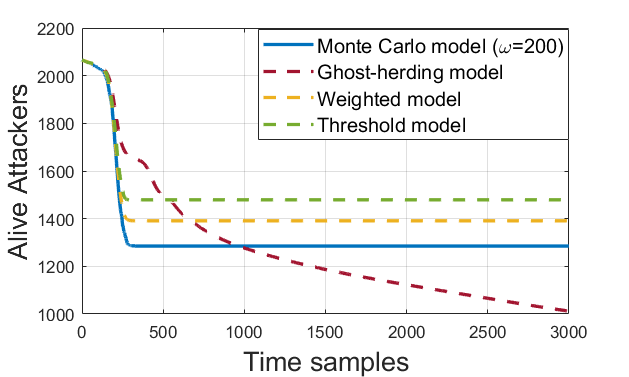}
\includegraphics[width=0.325\textwidth]{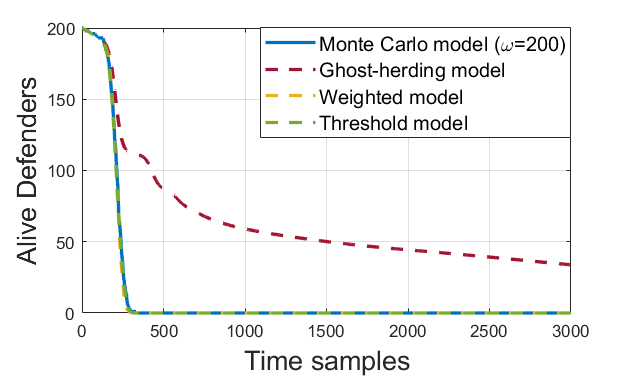}
\caption{Comparison of the performance of the 3 proposed models for optimization with the Monte Carlo simulation model in a scenario of 2066 attackers and 200 defenders}
\label{fig:compare-Ps}
\end{figure*}


Figure~\ref{fig:compare-Ps} shows results from each type of simulation, specifically the mean survival probabilities of attackers, defenders, and the HVU. This figure shows that $P_2$ and $P_3$ have similar results to Monte Carlo, which is typical of all simulations. However, $P_1$ does not agree with Monte Carlo, $P_2$, or $P_3$. Instead, $P_1$ greatly overestimates the probability of HVU survival, which is also typical of all simulations. Physically, this has an obvious explanation: $P_1$ had no coupling between attrition and spatial dynamics, so attackers would still try to avoid defenders who were in their path, even if these defenders had a very low probability of survival (ghost herding). In contrast, the modeling frameworks corresponding to $P_2$ and $P_3$ reduced the spatial interactions as survival probability decreased. So, even with the simple coupling between agent dynamics and agent attrition used in $P_2$ and $P_3$  these approximations agreed relatively well with the stochastic results produced by Monte Carlo simulation. This result highlights our main point in this paper, which is that modeling and control frameworks for adversarial autonomy must include attrition, and the attrition modeling should be coupled to the spatial dynamics of the agents.

\section{TRADE-OFF ANALYSIS}\label{sec:tradeoff}

\begin{figure*}
\centering
\includegraphics[trim=0mm 0mm 0mm 0mm, clip, width=0.50\textwidth]{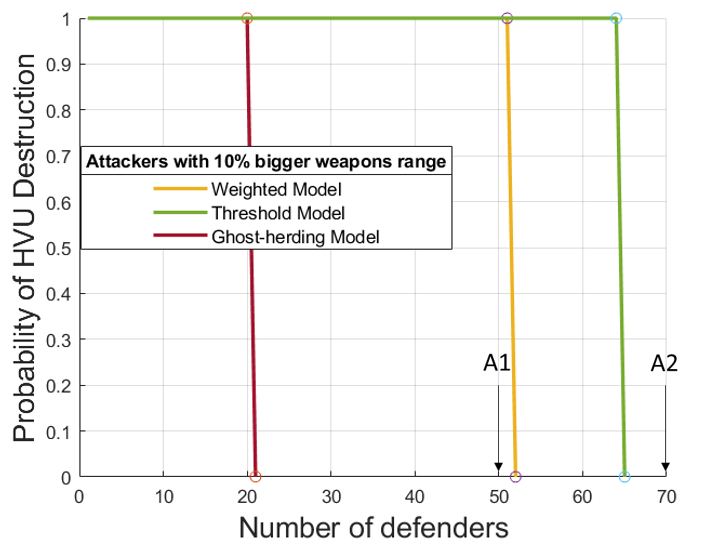}
\includegraphics[trim=0mm 0mm 0mm 0mm, clip,width=0.49\textwidth]{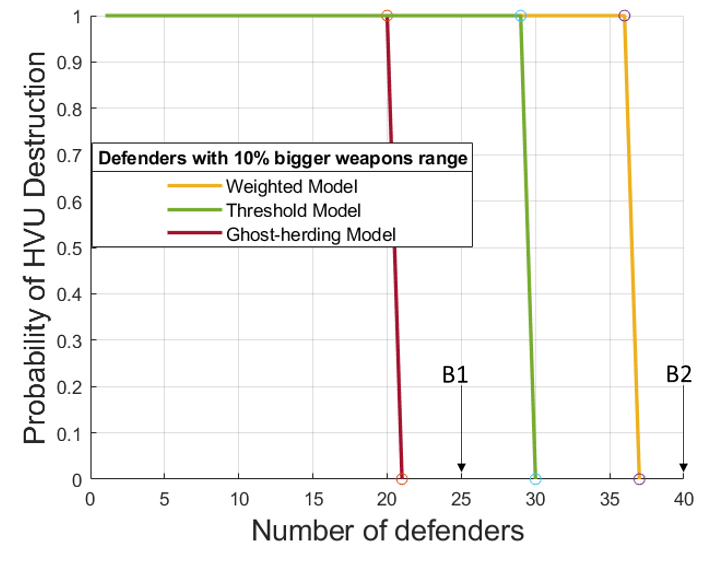}
\caption{Optimized cost versus number of defenders for all three proposed models for optimization against a swarm of 50 assets}
\label{fig:all-checkpoints}
\end{figure*}

In the previous section, we compared the three deterministic models with the stochastic Monte Carlo benchmark for a highly asymmetric engagement involving 2066 attackers and 200 defenders. That comparison established the principal modeling result: the coupled models $P2$ and $P3$ reproduce the stochastic outcome much more faithfully than the uncoupled model $P1$, which can exploit the unphysical ``ghost-herding'' mechanism. The purpose of this section is different. Here we use the models to examine a resource-allocation question that is directly relevant to defensive planning: 
\emph{how does the probability of HVU survival change as the number of defenders is increased, and how sensitive is the inferred defender requirement to the attrition model and to relative weapon capability?}

We consider a swarm of 50 attackers and vary the number of defenders $M$ from 1 to 70. For each value of $M$, the defender trajectories are re-optimized independently using each of the three deterministic models. The optimized objective, $J=1-Q_0(t_f)$, is the predicted probability of HVU destruction; equivalently, $Q_0(t_f)=1-J$ is the predicted HVU survival probability. Thus, Fig.~\ref{fig:all-checkpoints} can be interpreted as a family of \emph{resource--survival curves}: moving from left to right increases the defensive resources, while the sharp decrease in $J$ identifies the region in which the optimized defense changes from unsuccessful to successful. We refer to the smallest defender count for which the optimized model predicts essentially zero probability of HVU destruction as the \emph{critical defender count}. This quantity is model-dependent and is used here as a compact measure of the trade-off between defensive resources and mission success.

Two cases are examined to separate the effect of force size from the effect of relative weapon capability. In the A-type simulations [Fig.~\ref{fig:all-checkpoints}(a)], the attackers have a 10\% larger weapon range than the defenders. In the B-type simulations [Fig.~\ref{fig:all-checkpoints}(b)], the advantage is reversed and the defenders have a 10\% larger weapon range. All other aspects of the comparison are kept the same. Consequently, the horizontal displacement of the transition in Fig.~\ref{fig:all-checkpoints} directly measures how each model converts a modest change in weapon effectiveness into a change in the number of defenders required for HVU protection.

\subsection{Defender-number/HVU-survival trade-off}
For the A-type engagement, the uncoupled model $P1$ predicts a critical defender count of only 21. In contrast, the weighted model $P2$ requires approximately 52 defenders and the threshold model $P3$ approximately 65 defenders. The difference is operationally significant: the uncoupled model would suggest that a defending force substantially smaller than the 50-agent attacking swarm is sufficient even though the attackers also possess the weapon-range advantage. The coupled models instead predict that a force comparable to, or larger than, the attacking swarm is required. The discrepancy is not simply a numerical offset in the objective. It changes the resource level at which the model predicts that the mission becomes feasible.

The origin of this discrepancy is the ghost-herding mechanism. In $P1$, defenders continue to contribute to the attackers' avoidance dynamics even after their survival probabilities have become very small. The optimizer can therefore obtain protection by positioning defenders so that their repulsive influence redirects the swarm, without paying the physical consequence that destroyed defenders should cease to influence the attackers. Increasing the number of such defenders strengthens this artificial barrier. As a result, the apparent resource requirement in $P1$ is governed largely by the geometry of herding rather than by the ability of the defending force to remain effective under attrition. The critical count of 21 should therefore not be interpreted as a credible force-sizing prediction.

The difference between $P2$ and $P3$ in the A-type case is also informative. Both models remove the persistent influence of low-survival defenders, but they do so differently. In $P2$, an agent's dynamical influence decreases continuously with its survival probability; in $P3$, that influence is retained until the prescribed survival threshold is crossed and is then removed. Their predicted critical counts, 52 and 65, respectively, show that the particular approximation used to represent attrition can still affect quantitative force-sizing conclusions near the success/failure boundary. Thus, agreement that attrition must be coupled to the dynamics does not imply that all coupled approximations give identical resource requirements.

The B-type engagement provides an important sensitivity test. When the defenders are given a 10\% weapon-range advantage, the critical counts predicted by the coupled models decrease substantially: from 52 to 37 for $P2$ and from 65 to 30 for $P3$. This is the expected qualitative behavior. More effective defenders can engage attackers at greater range, reducing the number of defenders required to maintain a sufficiently effective defensive screen and protect the HVU. The change in critical count therefore demonstrates that the coupled models respond to a physically meaningful change in combat effectiveness.

By contrast, $P1$ again predicts a critical defender count of 21, exactly the same value obtained when the attackers possess the range advantage. This insensitivity is particularly revealing. Reversing a 10\% weapon-range advantage should alter the balance between attrition and protection, yet the resource requirement predicted by $P1$ does not change. The result confirms that the apparent success of the uncoupled model is dominated by ghost herding rather than by the modeled exchange of attrition. In this sense, the trade-off curves provide a stronger diagnostic than a single engagement: they show not only that $P1$ can give an incorrect outcome, but also that it can give a qualitatively incorrect \emph{sensitivity} of the required force size to weapon effectiveness.

The ordering of $P2$ and $P3$ also changes between the two cases. With stronger attackers, $P3$ is more conservative than $P2$ (65 versus 52 defenders), whereas with stronger defenders it is less conservative (30 versus 37 defenders). We therefore do not interpret either coupled approximation as uniformly conservative. Rather, the threshold operation in $P3$ changes the timing at which agents cease to influence the engagement, and that timing interacts with the relative attrition rates and optimized trajectories. This observation motivates the Monte Carlo comparisons below: the relevant question is not which deterministic model is always more pessimistic, but which model more faithfully reproduces the stochastic engagement over the range of force ratios of interest.

From a force-planning perspective, Fig.~\ref{fig:all-checkpoints} therefore contains two distinct trade-offs. The first is the expected trade-off between \emph{resources and survival}: increasing $M$ eventually drives the predicted probability of HVU destruction from near one to near zero. The second is a \emph{modeling trade-off}: the estimated location of that transition depends strongly on whether attrition is allowed to modify the spatial interaction dynamics. For the scenarios considered here, this modeling choice changes the inferred critical force size by tens of defenders. Hence a model that predicts plausible trajectories but treats destroyed agents as dynamically influential can lead to materially misleading resource-allocation conclusions.

\subsection{Pareto interpretation of the resource--survival curves}
The same results admit a useful multi-objective interpretation. For a fixed attrition model $P_r$, $r\in\{1,2,3\}$, let $\mathcal{U}_r(M)$ denote the feasible set of defender trajectories when $M$ defenders are available, and define the best predicted terminal HVU survival probability as
\begin{equation}
Q_{0,r}^{\star}(M)
\triangleq
\max_{\mathbf{u}\in\mathcal{U}_r(M)} Q_0(t_f;\mathbf{u},M),
\label{eq:QstarPareto}
\end{equation}
or, equivalently, the optimized destruction probability
\begin{equation}
J_r^{\star}(M)=1-Q_{0,r}^{\star}(M).
\label{eq:JstarPareto}
\end{equation}
Thus each optimization in Fig.~\ref{fig:all-checkpoints} produces a resource--performance pair
\begin{equation}
\big(M,J_r^{\star}(M)\big),
\end{equation}
where both coordinates are quantities to be minimized: fewer defenders are preferred, and a smaller destruction probability is preferred. A pair $\big(M_1,J_r^{\star}(M_1)\big)$ dominates $\big(M_2,J_r^{\star}(M_2)\big)$ if
\begin{equation}
M_1\le M_2,\qquad
J_r^{\star}(M_1)\le J_r^{\star}(M_2),
\end{equation}
with at least one strict inequality. The nondominated points obtained from the discrete sweep over $M$ therefore form an empirical, model-dependent Pareto frontier between defensive resources and predicted mission performance. We use the term \emph{Pareto-type resource--survival frontier} here because $M$ is swept externally rather than optimized simultaneously with the trajectories in a single bi-objective program.

This interpretation provides two complementary force-sizing questions. If a required HVU survival level $Q_{\rm req}$ is prescribed, the minimum force satisfying that requirement is
\begin{equation}
M_r^{\star}(Q_{\rm req})
=
\min\left\{M:\;Q_{0,r}^{\star}(M)\ge Q_{\rm req}\right\}.
\label{eq:MstarPareto}
\end{equation}
Conversely, if the available defender force is capped at $M_{\max}$, the relevant question is the maximum survival probability attainable with that budget. The ``critical defender count'' used above is a special case of \eqref{eq:MstarPareto} in which $Q_{\rm req}$ is chosen close to one. This makes Fig.~\ref{fig:all-checkpoints} more than a parameter sweep: it is a discrete approximation of the resource--performance frontier that a planner would use either to determine the minimum force required for a desired survival level or to determine the survival level achievable with a fixed force budget.

The Pareto viewpoint also sharpens the interpretation of ghost herding. Because $P1$ permits defenders with negligible survival probability to continue exerting repulsive influence on the attackers, it generates artificially favorable resource--performance pairs: high predicted HVU survival can be achieved with unrealistically small $M$. In geometric terms, the apparent $P1$ frontier is shifted toward the favorable low-resource/low-destruction region. The coupled models $P2$ and $P3$ move the inferred frontier toward larger defender counts because they remove, continuously or discontinuously, the dynamical influence of defenders as survival decreases. The modeling error in $P1$ therefore does not merely bias a single trajectory or terminal probability; it can distort the entire inferred Pareto trade-off and consequently lead to systematically optimistic force-allocation decisions.

It is important, however, not to interpret the curves as a formal continuous Pareto front for the underlying stochastic problem $P0$. The defender number is integer-valued, only the sampled values $M=1,\ldots,70$ are considered, and each point is obtained from a separate single-objective trajectory optimization. Moreover, $P1$--$P3$ are different approximations of $P0$. Accordingly, the curves are best viewed as \emph{model-dependent empirical Pareto frontiers}. Their principal value here is comparative: they reveal how the assumed treatment of attrition changes the resource--survival trade-off that would be presented to a decision maker.

\subsection{Checkpoint validation of the trade-off curves}
The optimized curves in Fig.~\ref{fig:all-checkpoints} summarize terminal mission outcome, but they do not show how the engagement evolves in time. We therefore examine four checkpoints, A1, A2, B1, and B2, chosen on opposite sides of the transition regions. A1 and B1 correspond to defender counts for which the coupled models predict insufficient protection, whereas A2 and B2 correspond to defender counts for which they predict successful protection. This comparison tests whether the force-sizing conclusions inferred from the terminal cost are also reflected in the time histories of HVU, attacker, and defender survival.

For each checkpoint, defender trajectories are first obtained using the weighted model $P2$. Those same trajectories and identical initial conditions are then used in $P1$, $P2$, $P3$, and an ensemble of 200 Monte Carlo realizations. Holding the defender trajectories fixed is important: differences among the resulting survival histories can then be attributed primarily to the attrition/dynamics model rather than to different optimized trajectories. The checkpoint study therefore complements the optimization curves by asking whether the deterministic models predict the stochastic outcome correctly when evaluated on the same defensive strategy.

\begin{figure*}
\centering
\begin{subfigure}{0.49\textwidth}
\centering\includegraphics[width=\textwidth]{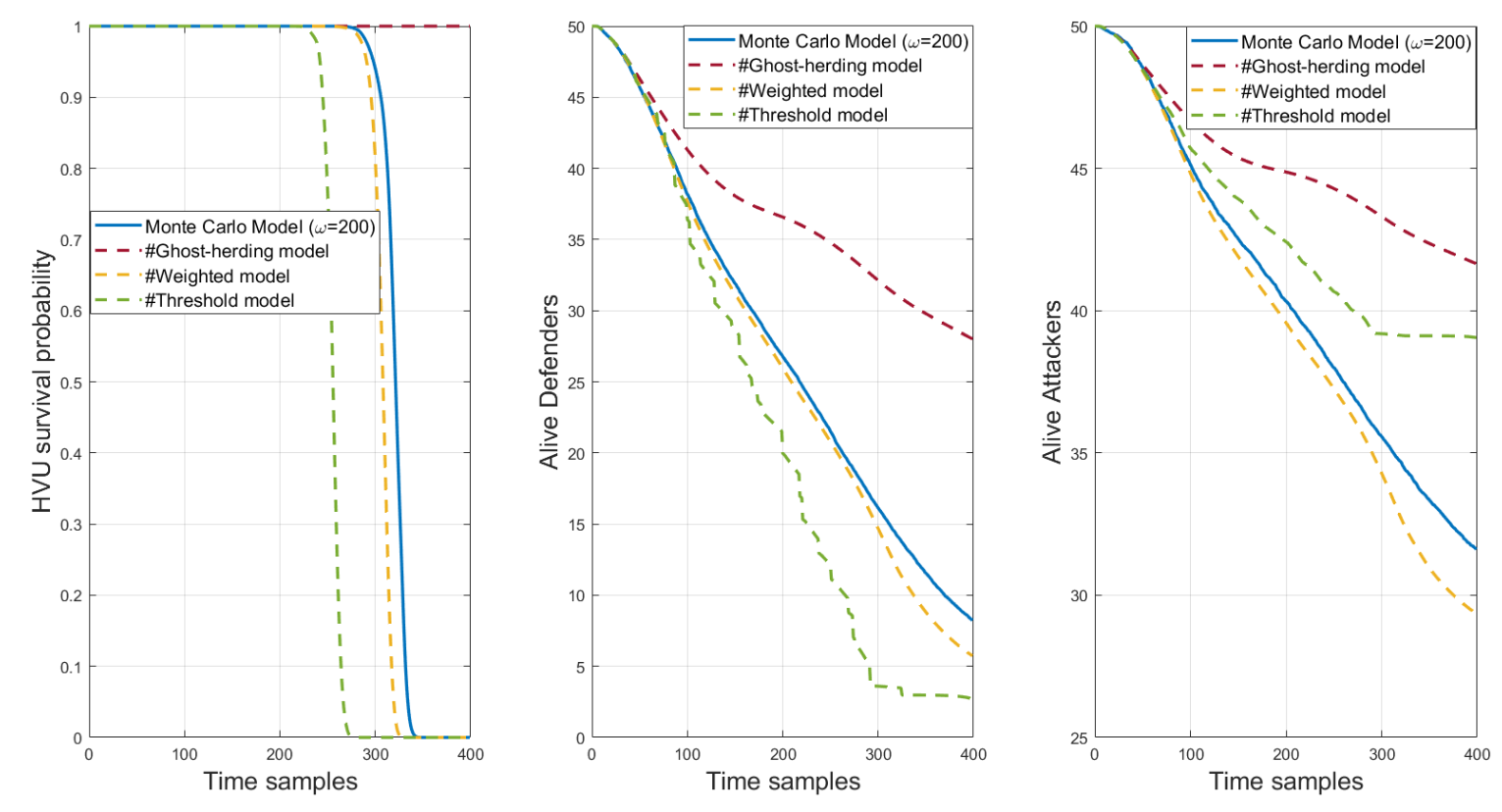}
\caption{A1: insufficient defenders.}\label{fig:checkpoint-A1}
\end{subfigure}\hfill
\begin{subfigure}{0.49\textwidth}
\centering\includegraphics[width=\textwidth]{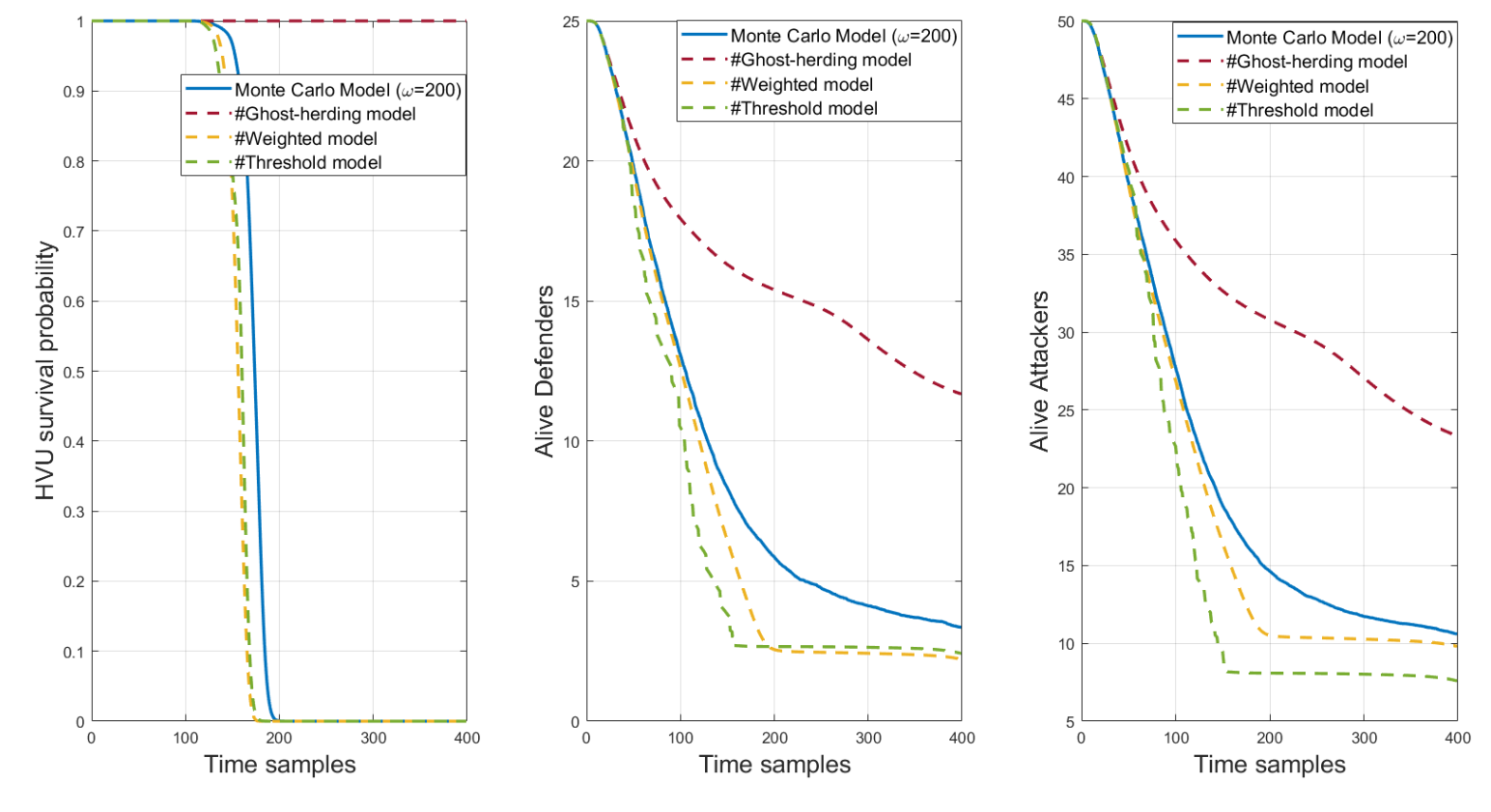}
\caption{B1: insufficient defenders.}\label{fig:checkpoint-B1}
\end{subfigure}
\caption{Checkpoint validation for cases in which the coupled models predict insufficient defensive resources.}
\end{figure*}

\begin{figure*}
\centering
\begin{subfigure}{0.49\textwidth}
\centering\includegraphics[width=\textwidth]{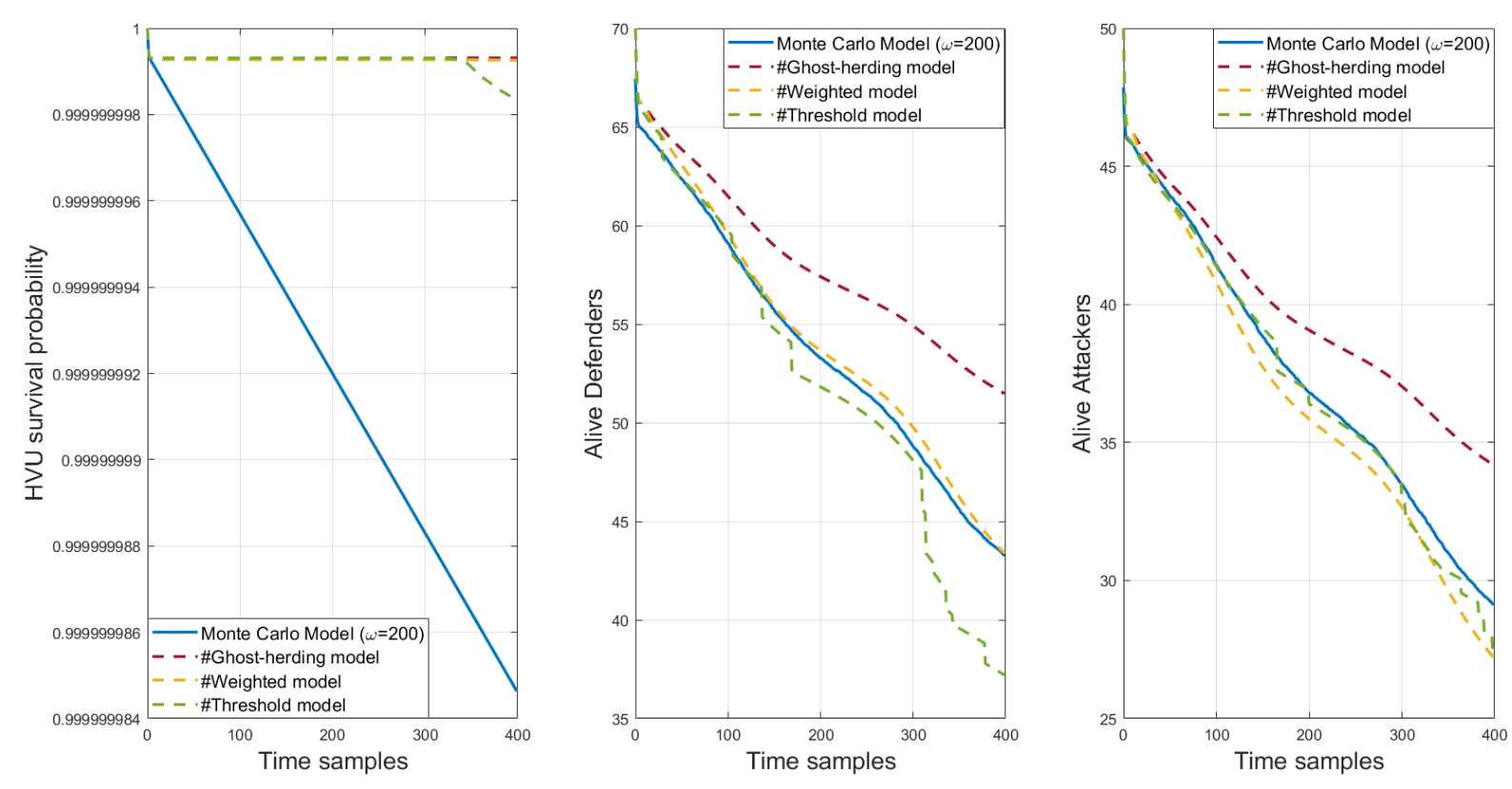}
\caption{A2: sufficient defenders.}\label{fig:checkpoint-A2}
\end{subfigure}\hfill
\begin{subfigure}{0.49\textwidth}
\centering\includegraphics[width=\textwidth]{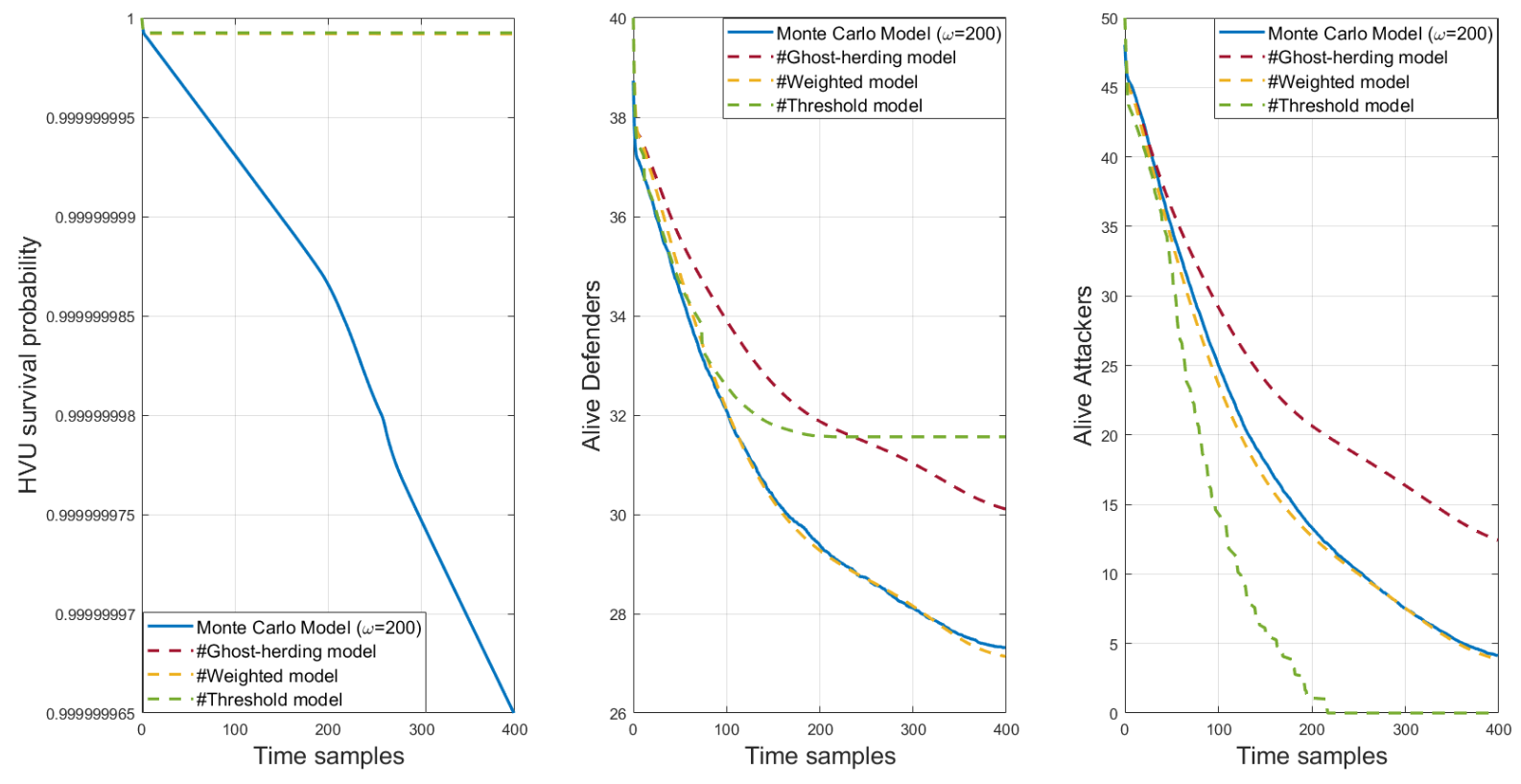}
\caption{B2: sufficient defenders.}\label{fig:checkpoint-B2}
\end{subfigure}
\caption{Checkpoint validation for cases in which the coupled models predict sufficient defensive resources.}
\end{figure*}

The first two checkpoints, A1 and B1, represent cases where the weighted and the threshold models $P2$ and $P3$ both predict that the HVU will be destroyed. The checkpoints A2 and B2 represent cases where weighted and threshold models both predict that the number of defenders is sufficient to protect the HVU. We also compare with the ensemble average of 200 Monte Carlo simulations. We note two important features of the data shown in Figs.~\ref{fig:checkpoint-A1} through \ref{fig:checkpoint-B2}. First, the HVU survival probability versus time shows good agreement between the two coupled models (weighted and threshold, $P2$ and $P3$) and the Monte Carlo simulations. In the two cases (A1 and B1) where the HVU is destroyed, it happens at roughly the same time for all three cases. The uncoupled model $P1$ always predicts HVU survival, due to the ``ghost-herding'' problem. Second, for the survival probability of the attackers and defenders versus time, the weighted model shows the closest agreement with the threshold model. This agreement is particularly striking for checkpoint B2, shown in Fig.~\ref{fig:checkpoint-B2}, but it is a generic feature to all four checkpoints shown.

\subsection{Scope of validation and role of physical experiments}
The purpose of the present validation is to assess the treatment of stochastic attrition within the optimal-control model, rather than to validate a particular air vehicle, sensor suite, communication architecture, or weapon system. For this reason, the stochastic formulation $P0$ is the appropriate reference model for evaluating the deterministic approximations $P1$--$P3$. The Monte Carlo study provides a stochastic benchmark: trajectories generated by the deterministic optimization are held fixed and evaluated over repeated realizations in which agents are actually removed from the active population according to the attrition model. Thus, the comparison directly tests the question posed in this paper---whether a computationally tractable deterministic model preserves the consequences of stochastic agent removal.

Physical swarm experiments address a complementary validation question. They are particularly valuable for establishing the fidelity of vehicle dynamics, sensing, communication, collision avoidance, and inter-agent interaction mechanisms. Indeed, experimental studies have already demonstrated the physical feasibility of aerial swarm-defense and herding concepts; for example, the StringNet framework was validated experimentally using quadrotors in \cite{chipade2021stringnet}. Such experiments support the feasibility of the underlying defender--attacker interaction mechanisms, but they do not by themselves validate the probabilistic attrition law used here. A laboratory experiment that emulates attrition by probabilistically declaring vehicles inactive would still obtain those removal events from a software-specified stochastic model and would therefore amount to a hardware-in-the-loop realization of the same attrition assumption.

There is also an important issue of scale. The methodology is intended for engagements involving populations that can reach hundreds or thousands of agents, as illustrated by the 2066-attacker/200-defender case considered above. Full physical replication at this scale is presently impractical. We therefore regard Monte Carlo validation against $P0$ as the appropriate validation of the \emph{modeling approximation} studied in this paper, while experimental calibration of the motion, interaction, and attrition laws remains a separate and important step before the framework is used to make quantitative predictions for a specific physical system. In particular, the critical defender counts and empirical Pareto frontiers reported in this section should be interpreted as \emph{model-based resource--survival trade-offs}, not as experimentally calibrated force requirements for a particular real-world engagement. Future work should combine empirically calibrated interaction and attrition models with hardware-in-the-loop and scaled multi-vehicle experiments.



\section{Conclusions}\label{sec:conclusions}

In this paper we have addressed the question of optimal motion planning for large-scale, coordinated autonomous systems in adversarial conditions, by which we mean that it is likely that some agents will be lost while trying to accomplish the task. We studied this scenario using the framework of optimal control. We specifically studied the scenario of a multi-agent team defending a high-value unit against a swarm attack. 

We observed several key challenges associated with using optimal control to study such a large-scale motion planning problem in adversarial conditions. The first challenge revolves around a modeling framework to study the scenario. We demonstrate a coupled framework between dynamics, Eqs.~\eqref{eqn:attackers} and \eqref{eqn:defenders}, and attrition modeling, Eqs.~\eqref{eqn:Q(t)-all}, is an appropriate and highly flexible way to model such engagements.

The second challenge involves formulation on an optimal control problem, which in principle could be used to find the best set of trajectories to accomplish a given task. We proposed a novel optimal control problem, $P0$, that explicitly includes random reduction of an index set of surviving agents in time. 

However, since there is no known solution to such a problem, the third challenge involved asking whether there are approximate problems that could be proposed and solved that would faithfully mimic the stochastic optimal control problem. To this end, we proposed three {\color{black} non-random} approximations that can be solved using direct methods of optimal control. By considering a case study of defending agents protecting an HVU from an attacking swarm, we showed that these approximations can be solved and that they give results that are consistent with the stochastic problem, especially if the attrition and spatial dynamics are coupled. 

We note that our results assume that the coordinating strategy and weapons capabilities of the attacking swarm are known or can be estimated. Recent work by Peltier et al.~\cite{peltier2025swarmclassification,peltier2025robustclassification} demonstrates that attacking-swarm characteristics and tactics can be inferred from observed trajectories using neural-network time-series classification, including under variations in defender number, defender motion, and measurement noise. The latter work further optimizes defender trajectories to elicit adversarial responses that improve tactic classification. These capabilities are complementary to the framework developed here: classification provides information about how the adversary is coordinating, whereas the present framework determines how the engagement dynamics should evolve as agents are removed through attrition. This suggests a natural future architecture in which online tactic classification is coupled with attrition-aware trajectory optimization, allowing the defensive strategy to adapt jointly to the inferred adversarial tactic and the evolving surviving-agent population. Parameter uncertainty can also be incorporated using methods such as in \cite{walton_IJC}. The framework we describe here can be applied to an entire class of adversarial autonomy problems. For example, attrition could result from many factors, including environmental or terrain-related causes. To the best of the authors' knowledge, prior adversarial-swarm defense studies have not explicitly identified or analyzed the inconsistency that occurs when attrition is modeled while destroyed agents continue to exert spatial influence. The ghost-herding phenomenon identified here shows that this inconsistency can distort optimized trajectories, predicted HVU survival, and the inferred resource--survival frontier. Thus, the contribution is not simply another swarm-defense strategy, but a modeling framework for ensuring that stochastic attrition is consistently reflected in the dynamics used for optimal defensive decision making.

The numerical results in this paper validate the deterministic approximations against the stochastic model $P0$; they do not constitute experimental calibration of a particular physical engagement. Consequently, the defender counts and resource--survival frontiers should be interpreted conditionally on the assumed motion, interaction, and attrition models. Before quantitative force-sizing conclusions are transferred to a specific system, those component models should be calibrated using appropriate experimental or operational data. Future work should therefore include hardware-in-the-loop and scaled multi-vehicle experiments, together with improved approximation methods ($P1$, $P2$, and $P3$) and stochastic optimization techniques that may ultimately permit $P0$ to be treated directly.

The numerical results in this paper validate the deterministic approximations against the stochastic model $P0$; they do not constitute experimental calibration of a particular physical engagement. Consequently, the defender counts and resource--survival frontiers should be interpreted conditionally on the assumed motion, interaction, and attrition models. Before quantitative force-sizing conclusions are transferred to a specific system, those component models should be calibrated using appropriate experimental or operational data. Future work should therefore include hardware-in-the-loop and scaled multi-vehicle experiments, together with improved approximation methods ($P1$, $P2$, and $P3$) and stochastic optimization techniques that may ultimately permit $P0$ to be treated directly.

\section*{Acknowledgment}
This work was supported in part by the ONR Science of Autonomy program under the grant N0001425GI01578 and by NPS CRUSER program. 
\ifCLASSOPTIONcaptionsoff
  \newpage
\fi

\bibliographystyle{IEEEtran}
\bibliography{references}  
\end{document}